\documentclass[a4paper]{article}
\usepackage{amsfonts}
\usepackage{amssymb}
\usepackage{amsmath}
\usepackage{latexsym}
\usepackage{graphicx}
\usepackage{calrsfs}
 \newtheorem{thm}{Theorem}[section]

 \newtheorem{rem}[thm]{Remark}

\begin{document}

\parbox{1mm}

\begin{center}
{\bf {\sc \Large On Toeplitz localization operators}}
\end{center}

\vskip 12pt

\begin{center} {\bf Ondrej HUTN\'IK and M\'aria HUTN\'IKOV\'A}
\end{center}

\vskip 24pt

\hspace{5mm}\parbox[t]{10cm}{\fontsize{9pt}{0.1in}\selectfont\noindent{\bf
Abstract.} We present a unified approach to study properties of
Toeplitz localization operators based on the Calder\'on and Gabor
reproducing formula. We show that these operators with functional
symbols on a plane domain may be viewed as certain
pseudodifferential operators (with symbols on a line, or certain
compound symbols). } \vskip 24pt

\section{Introduction and preliminaries}

A starting point for the construction of time-frequency
localization or filter operators are the famous reproducing
formulas of Calder\'on (in wavelet analysis) and of Gabor (in
time-frequency analysis). In this paper we will work with both the
reproducing formulas and therefore we introduce the following
unified notation. Write $$f=\int_{G} \langle f,\Psi_{\zeta}\rangle
\Psi_{\zeta} \,\mathrm{d}\zeta,$$ where $G$ denotes either the
half-plane $\mathbb{R}^{2}_+ = \mathbb{R}_+\times \mathbb{R}$, or
the whole plane $\mathbb{R}^{2}$. In the first (wavelet) case
$\zeta=(u,v)$, $u>0$, $v\in\mathbb{R}$, and
$$\psi_{u,v}(x)=\frac{1}{\sqrt{u}}\, \psi\left(\frac{x-v}{u}\right)$$ is the
action of the group $\mathbb{R}^{2}_+$ on $L_2(\mathbb{R})$, where
$\mathbb{R}^2_+$ is equipped with the hyperbolic measure
$\mathrm{d}\zeta = u^{-2}\,\mathrm{d}u\mathrm{d}v$. Here (the
real-valued) \textit{admissible wavelet} is the function $\psi\in
L_2(\mathbb{R})$ satisfying the condition
$$\int_{\mathbb{R}_{+}} |\hat{\psi}(t\xi)|^{2}\frac{\mathrm{d}t}{t}=1$$ for
almost every $\xi\in\mathbb{R}$, where $\hat{\psi}$ stands for the
Fourier transform (in the time-frequency convention) $\mathcal{F}:
L_{2}(\mathbb{R})\to L_{2}(\mathbb{R})$ given by
$$\mathcal{F}\{f\}(\xi) = \hat{f}(\xi)=\int_{\mathbb{R}}
f(x)\,\mathrm{e}^{-2\pi \mathrm{i}x\xi}\,\mathrm{d}x.$$ In the
second (time-frequency) case $\zeta=(q,p)$, $q,p\in\mathbb{R}$
with
$$\phi_{q,p}(x) = \mathrm{e}^{2\pi \mathrm{i}px} \phi(x-q)$$ being the action of
$\mathbb{R}^{2}$ on $L_2(\mathbb{R})$, and
$\mathrm{d}\zeta=\mathrm{d}q\mathrm{d}p$ being the measure on
$\mathbb{R}^2$. The \textit{admissible window} is the function
$\phi\in L_2(\mathbb{R})$ satisfying
$\|\phi\|_{L_2(\mathbb{R})}=1$. In what follows the symbol $\Psi$
always means either an admissible wavelet $\psi\in
L_2(\mathbb{R})$, or an admissible window $\phi\in
L_2(\mathbb{R})$.

If $L_2(G,\mathrm{d}\zeta)$ denotes the Hilbert space of all
square-integrable complex-valued functions on $G$, then for a
fixed $\Psi \in L_{2}(\mathbb{R})$ the functions $W_{\Psi}f$ on
$G$ of the form
$$(W_{\Psi} f)(\zeta)=\langle f, \Psi_{\zeta}\rangle, \quad f\in
L_{2}(\mathbb{R}),$$ form a reproducing kernel Hilbert space
$W_{\Psi}(L_{2}(\mathbb{R}))$. Then the transform $W_{\Psi}:
L_2(\mathbb{R}) \to L_{2}(G, \mathrm{d}\zeta)$ is an isometry, and
the integral operator $P_{\Psi}: L_{2}(G, \mathrm{d}\zeta) \to
L_{2}(G, \mathrm{d}\zeta)$ given by
$$(P_{\Psi}F)(\eta)=\int_{G} F(\zeta) \langle \Psi_\eta, \Psi_\zeta
\rangle\,\mathrm{d}\zeta, \quad F\in L_{2}(G, \mathrm{d}\zeta),$$
is the orthogonal projection onto $W_{\Psi}(L_{2}(\mathbb{R}))$,
where $\langle \Psi_{\eta}, \Psi_{\zeta}\rangle$ is the
reproducing kernel in $W_{\Psi}(L_{2}(\mathbb{R}))$. For a given
bounded function $a$ on $G$ define the \textit{Toeplitz
localization operator} $T_a^{\Psi}$ with symbol $a$ as follows
$$T_a^\Psi: f\in W_\Psi(L_2(\mathbb{R})) \longmapsto P_\Psi(af)\in
W_\Psi(L_2(\mathbb{R})).$$ In wavelet case the Toeplitz operator
$T_a^\psi$ is usually called the \textit{Calder\'on-Toeplitz
operator}, whereas in the case of time-frequency analysis the
operator $T_a^\phi$ is called the \textit{Gabor-Toeplitz
operator}.

In this paper we underline some interesting features of spaces of
transforms $W_\Psi(L_2(\mathbb{R}))$ and Toeplitz localization
operators acting on them. In fact, with the above notation we
provide a unified approach to both cases and give a more natural
construction of unitary operators which does not use the
decomposition of $L_2(G, \mathrm{d}\zeta)$ onto spaces
$W_\Psi(L_2(\mathbb{R}))$ as it was done in~\cite{hutnik}
and~\cite{huthut}. Indeed, in Section~\ref{section2} according to
the general scheme presented in~\cite{Q-BV} we give the
construction of the unitary operator $R_\Psi$ which is an exact
analog of the Bargmann transform mapping the Fock space
$F_2(\mathbb{C}^n)$ of Gaussian square-integrable entire functions
on $\mathbb{C}^n$ onto $L_2(\mathbb{R}^n)$, see~\cite{bargmann}.
Then, via $R_\Psi$, the Toeplitz localization operators
$T_a^{\Psi}: W_\Psi(L_2(\mathbb{R}))\to W_\Psi(L_2(\mathbb{R}))$
can be identified with certain pseudodifferential operators
$$\mathfrak{C}_a^\Psi: = R_\Psi T_a^{\Psi} R_\Psi^*:
L_2(\mathbb{R})\to L_2(\mathbb{R}).$$ This passing from
$T_a^{\Psi}$ to $\mathfrak{C}_a^{\Psi}$ is nothing but an analog
of the Berezin reducing of (Toeplitz) operators with anti-Wick
symbols on the Fock space $F_2(\mathbb{C}^n)$ to Weyl
pseudodifferential operators on $L_2(\mathbb{R}^n)$,
see~\cite{berezin} for further details.

In particular, in Section~\ref{section3.1} the above mentioned
observation is applied to operator symbols $a(r,s): G \to
\mathbb{C}$ that are only depending on the variable $r$. In this
case (cf. Theorem~\ref{CTO1}) the operator $\mathfrak{C}_a^\Psi$
is simply a multiplication operator with explicitly computable
symbol $\gamma_a^\Psi$. As a consequence the boundedness (also for
the case where $a$ is unbounded!) and the spectrum of $T_a^\Psi$
is precisely characterized in terms of function $\gamma_a^\Psi$.
Moreover, for a fixed $\Psi\in L_2(\mathbb{R})$ the space of
Toeplitz localization operators $T_a^\Psi$ with bounded symbols
$a$ depending only on $r$ generates a commutative
$C^*$-sub-algebra of the Toeplitz $C^*$-algebra. In
Theorem~\ref{thm3.6} the commutative $C^*$-algebra generated by
operators $T_a^\Psi$ with symbols $a(r,s) = \alpha(r)$ and
$\alpha$ in a certain space of piecewise constant functions is
shown to be isometrically isomorphic to an explicitly given
algebra of continuous functions.

In Section~\ref{section3.2} the case of an operator symbol $a(r,
s) = \beta(s)$ depending only on the second variable $s$ of $G$ is
studied. In this case $\mathfrak{C}_a^\Psi$ has the form of a
certain integral operator on $L_2(\mathbb{R})$. Finally,
Section~\ref{section3.3} treats the mixed case $a(r, s) =
\alpha(r)\beta(s)$ which leads to pseudodifferential operators
$\mathfrak{C}_a^\Psi$ with an explicitly computable compound
(double) symbol. These results provide an interesting tool for
further study of Toeplitz localization operators via investigating
their unitary equivalent images in the classes of
pseudodifferential operators.

\section{Bargmann-type transform}\label{section2}

In what follows let $\Psi \in L_2(\mathbb{R})$ be fixed. In order
to construct unitary operators which will be used to study
Toeplitz localization operators, we represent the Hilbert space
$L_2(G,\mathrm{d}\zeta)$ as a tensor product in the form
$$L_2(G, \mathrm{d}\zeta) =  L_{2}(G_1, \mathrm{d}\zeta_1)\otimes L_{2}(G_2, \mathrm{d}\zeta_2),$$ where
$G_1=\mathbb{R}_+$, $G_2=\mathbb{R}$ with
$\mathrm{d}\zeta_1=u^{-2}\,\mathrm{d}u$,
$\mathrm{d}\zeta_2=\mathrm{d}v$ in the first (wavelet) case, and
$G_1=G_2=\mathbb{R}$ with $\mathrm{d}\zeta_1=\mathrm{d}q$,
$\mathrm{d}\zeta_2=\mathrm{d}p$ in the second (time-frequency)
case, respectively. Introduce the unitary operator
$$U_\Psi: L_2(G, \mathrm{d}\zeta) = L_{2}(G_1, \mathrm{d}\zeta_1)\otimes L_{2}(G_2, \mathrm{d}\zeta_2) \to
L_{2}(G_1, \mathrm{d}\zeta_1) \otimes L_{2}(G_2,
\mathrm{d}\zeta_2)$$ given by $U_\Psi = (I\otimes \mathcal{F}^{\pm
1})$, where the Fourier transform $\mathcal{F}=\mathcal{F}^{+1}$
corresponds to wavelet case, and the inverse Fourier transform
$\mathcal{F}^{-1}$ corresponds to time-frequency case. The image
$\Delta_\Psi$ of the space $W_{\Psi}(L_{2}(\mathbb{R}))$ under the
mapping $U_\Psi$ consists of all functions $F(z,\omega) =
f(\omega)\ell_\Psi(z,\omega)$, where $f \in L_2(\mathbb{R})$ and
$$\ell_\psi(u,\omega) =
\sqrt{u}\,\overline{\hat{\psi}(u\omega)}\quad \textrm{and}\quad
\ell_\phi(q,\omega) = \overline{\phi(\omega-q)},$$ respectively.
Clearly, for each $\omega\in G_2$ holds $$\ell_\Psi(\cdot,
\omega)\in L_2(G_1, \mathrm{d}\zeta_1)\,\,\,\,
\textrm{with}\,\,\,\, \|\ell_\Psi(\cdot, \omega)\|_{L_2(G_1,\,
\mathrm{d}\zeta_1)} = 1,$$ and thus we obviously have
$$\|F(z,\omega)\|_{\Delta_\Psi} =
\|f(\omega)\|_{L_{2}(G_2,\,\mathrm{d}\zeta_2)}.$$ Then the
operator $\Lambda_\Psi: L_{2}(G, \mathrm{d}\zeta) \to \Delta_\Psi$
given by $\Lambda_\Psi = U_\Psi P_{\Psi} U_\Psi^{*}$ has the
explicit form
$$(\Lambda_\Psi F)(z,\omega) = \ell_\Psi(z,\omega) \int_{G_1} F(t,\omega)\,
\overline{\ell_\Psi(t,\omega)}\,\mathrm{d}\zeta_1(t).$$ Thus,
$\textrm{Im}\,\Lambda_\Psi = \Delta_\Psi$. Moreover,
$\Lambda_\Psi^2 = \Lambda_\Psi$, and $\Lambda_\Psi$ is obviously
self-adjoint. Introduce the isometric imbedding $Q_\Psi: L_2(G_2,
\mathrm{d}\zeta_2) \to \Delta_\Psi$ by the rule $$(Q_\Psi
f)(z,\omega) = f(\omega)\ell_\Psi(z,\omega).$$ Then the adjoint
operator $Q_\Psi^*: L_2(G,\mathrm{d}\zeta)\to
L_2(G_2,\mathrm{d}\zeta_2)$ is given by
$$(Q_\Psi^* F)(\xi) = \int_{G_1} F(t,\xi)\,
\overline{\ell_\Psi(t,\xi)}\,\mathrm{d}\zeta_1(t),$$ and it is
easy to verify that the operators $Q_\Psi$ and $Q_\Psi^*$ provide
the following decomposition of identity on $L_2(G_2,
\mathrm{d}\zeta_2)$ and of orthogonal projection $\Lambda_\Psi$,
i.e.,
\begin{align*}
Q_\Psi^* Q_\Psi & = I_{\phantom{\Psi}}: L_2(G_2,
\mathrm{d}\zeta_2)\to L_2(G_2, \mathrm{d}\zeta_2), \\ Q_\Psi
Q_\Psi^* & = \Lambda_\Psi: L_2(G, \mathrm{d}\zeta) \to
\Delta_\Psi.\end{align*}The whole situation of constructed
operators is described on Figure~\ref{fig_operators}.

\begin{thm}\label{thmR_Psi}
The operator $R_\Psi=Q_\Psi^*U_\Psi$ maps the space $L_2(G,
\mathrm{d}\zeta)$ onto $L_2(G_2, \mathrm{d}\zeta_2)$, and the
restriction
$$R_\Psi\mid_{W_\Psi(L_2(\mathbb{R}))}: W_\Psi(L_2(\mathbb{R})) \to
L_2(G_2, \mathrm{d}\zeta_2)$$ is an isometrical isomorphism. The
adjoint
$$R_\Psi^*=U_\Psi^* Q_\Psi: L_2(G_2, \mathrm{d}\zeta_2)\to
W_\Psi(L_2(\mathbb{R}))\subset L_2(G, \mathrm{d}\zeta)$$ is an
isometrical isomorphism of $L_2(G_2, \mathrm{d}\zeta_2)$ onto the
subspace $W_\Psi(L_2(\mathbb{R}))$ of the space $L_2(G,
\mathrm{d}\zeta)$. Moreover,
\begin{align*}
R_\Psi R_\Psi^* & = I_{\phantom{\Psi}}: L_{2}(G_2,
\mathrm{d}\zeta_2) \to L_{2}(G_2, \mathrm{d}\zeta_2), \\ R_\Psi^*
R_\Psi & = P_\Psi: L_{2}(G, \mathrm{d}\zeta) \to
W_\Psi(L_2(\mathbb{R})).\end{align*}
\end{thm}

\begin{figure}
\begin{center}
\includegraphics[scale=0.9]{figoperators.1}
\end{center}\caption{Relationships among the constructed unitary operators}\label{fig_operators}
\end{figure}

In what follows we show that the Bargmann-type transform $R_\Psi$
essentially simplifies the previous computations made
in~\cite{hutnik} and~\cite{huthut}, and enables to obtain many
interesting results for the Toeplitz localization operators which
both cases share in common in a more transparent way.

\section{Toeplitz localization operators}

For each function $a(r,s)\in L_\infty(G, \mathrm{d}\zeta)$
consider the Toeplitz localization operator (TLO, for short)
$$T_a^{\Psi}: f\in W_\Psi(L_2(\mathbb{R}))\longmapsto P_\Psi(af)\in
W_\Psi(L_2(\mathbb{R})).$$ In our original
publications~\cite{hutnik} and~\cite{huthut} the operators
$R_\Psi$ and $R_\Psi^*$ were defined after applying the second
unitary operator, let us say $$V_\Psi: L_{2}(G_1,
\mathrm{d}\zeta_1)\otimes L_{2}(G_2, \mathrm{d}\zeta_2) \to
L_{2}(G_1, \mathrm{d}x) \otimes L_{2}(G_2, \mathrm{d}y)$$ given by
$$V_\psi: F(u,\omega) \mapsto
\frac{\sqrt{|y|}}{x}F\left(\frac{x}{|y|},y\right)\quad
\textrm{and}\quad V_\phi: F(q,\omega)\mapsto F(y-x,y).$$ Under the
operator $V_\Psi U_\Psi: L_{2}(G, \mathrm{d}\zeta)\to L_{2}(G_1,
\mathrm{d}x) \otimes L_{2}(G_2, \mathrm{d}y)$ we have obtained the
structural result saying "how much space occupies the subspace
$W_{\Psi}(L_{2}(\mathbb{R}))$ inside $L_2(G, \mathrm{d}\zeta)$",
see~\cite[Theorem 3.3]{hutnik} and~\cite[Theorem 1]{huthut} for
more details. Now the trick is that in comparison with our
previous approach the second operator $V_\Psi$ in both cases
\textit{is not needed} to study the TLO's $T_a^\Psi$, thus
providing a much easier way to the properties of $T_a^\Psi$. Of
course, the previous approach has its own advantages in connection
with understanding the structure of $W_{\Psi}(L_{2}(\mathbb{R}))$
inside $L_2(G, \mathrm{d}\zeta)$, as well as with study of certain
algebras of operators, e.g., algebras generated by operators of
the form $A_\Psi = aI+b P_\Psi$ acting on $L_2(G,
\mathrm{d}\zeta)$ with $a, b$ the bounded functions on $G$
depending only on the first coordinate. In what follows we
gradually apply the Bargmann-type transform $R_\Psi$ to the TLO
$T_a^\Psi$ with a symbol as a function of individual coordinates
of $G$.

\subsection{TLO's with symbols depending on the first
variable}\label{section3.1}

In what follows we consider the case where the symbol of the TLO
depends only on horizontal variable of $G$. This case is very
important because it gives rise to commutative operator algebras
with some interesting features.

\begin{thm}\label{CTO1}
If a measurable function $a(r,s)=\alpha(r)$ on $G$ does not depend
on $s$, then the TLO $T_{\alpha}^\Psi$ acting on
$W_{\Psi}(L_{2}(\mathbb{R}))$ is unitarily equivalent to the
multiplication operator $\gamma^\Psi_\alpha I$ acting on $L_2(G_2,
\mathrm{d}\zeta_2)$. The function $\gamma^\Psi_{\alpha}: G_2\to
\mathbb{C}$ is given by
$$\gamma^\Psi_{\alpha}(\xi)=\int_{G_1} \alpha(r)
|\ell_\Psi(r,\xi)|^2\,\mathrm{d}\zeta_1(r), \quad \xi\in G_2.$$
\end{thm}

\begin{proof} The operator $T_\alpha^\Psi$ is obviously unitarily
equivalent to the following operator
\begin{align*} R_\Psi T_{\alpha}^\Psi
R_\Psi^{*} & = R_\Psi P_{\Psi} \alpha(r) P_{\Psi} R_\Psi^{*} =
R_\Psi(R_\Psi^{*}R_\Psi)\alpha(r)(R_\Psi^{*}R_\Psi)R_\Psi^{*}
\\ & = (R_\Psi R_\Psi^{*})R_\Psi \alpha(r)R_\Psi^{*}(R_\Psi R_\Psi^{*}) =
R_\Psi \alpha(r) R_\Psi^{*} \\ & =
Q_\Psi^{*}\alpha(r)Q_\Psi,\end{align*}where the result of
Theorem~\ref{thmR_Psi} has been used. Finally, for $f\in L_2(G_2,
\mathrm{d}\zeta_2)$ we have
$$\left(Q_\Psi^{*}\alpha(r)Q_\Psi f\right)(\xi) = \int_{G_1}
\alpha(r) f(\xi) |\ell_\Psi(r,\xi)|^2\,\mathrm{d}\zeta_1(r) =
f(\xi)\cdot\gamma^\Psi_{a}(\xi), \quad \xi\in G_2,$$ where
$\gamma^\Psi_{\alpha}(\xi)=\int_{G_1} \alpha(r)
|\ell_\Psi(r,\xi)|^2\,\mathrm{d}\zeta_1(r)$ for $\xi\in G_2$.
\end{proof}

\begin{rem}\rm
The explicit form of the corresponding function $\gamma_a^\Psi$
for both cases is as follows $$ \gamma_\alpha^\psi(\xi) =
\int_{\mathbb{R}_+} \alpha(u)
|\hat{\psi}(u\xi)|^2\,\frac{\mathrm{d}u}{u}, \qquad
\gamma_\alpha^\phi(\xi) = \int_{\mathbb{R}} \alpha(q)
|\phi(\xi-q)|^2\,\mathrm{d}q, $$ where $\xi\in\mathbb{R}$. In
fact, the functions $\gamma^\psi_{\alpha}$ and
$\gamma_\alpha^\phi$ are constructed by putting a multiplier in
admissibility conditions for the wavelet $\psi$ and the window
$\phi$, respectively.
\end{rem}

Clearly, the result of Theorem~\ref{CTO1} opens an easy and direct
way to properties of TLO's with symbols depending only on first
variable. Since $T_\alpha^{\Psi}$ is unitarily equivalent to a
multiplication operator, then it is never compact. If
$a(r,s)=\alpha(r)$ is a bounded symbol, then the operator
$T_\alpha^{\Psi}$ is obviously bounded on
$W_\Psi(L_2(\mathbb{R}))$, and for its operator norm holds
$$\|T_\alpha^{\Psi}\| \leq \textrm{ess-sup}\,|\alpha(r)|.$$ Therefore, the
spaces $W_\Psi(L_2(\mathbb{R}))$ are natural and appropriate for
TLO's with bounded symbols. However, we may observe that the
result of Theorem~\ref{CTO1} suggests considering not only bounded
symbols, but also unbounded ones. In this case we obviously have

\begin{thm}\label{corgamma}
Let $a(r,s)=\alpha(r)$ be a measurable symbol on $G$. Then the TLO
$T_{\alpha}^\Psi$ is bounded on $W_{\Psi}(L_{2}(\mathbb{R}))$ if
and only if the corresponding function $\gamma_{\alpha}^\Psi(\xi)$
is bounded on $G_2$, and
$$\|T_\alpha^\Psi\| = \sup_{\xi\in G_2}
|\gamma_\alpha^\Psi(\xi)|.$$
\end{thm}

In this case we may also easily describe the spectrum of TLO as
follows.

\begin{thm}
The spectrum of a bounded TLO $T_\alpha^\Psi$ acting on
$W_\Psi(L_2(\mathbb{R}))$ with a measurable symbol
$a(r,s)=\alpha(r)$ coincides with its essential spectrum, and is
given by
$$\mathrm{sp}\,T_\alpha^\Psi =
\mathrm{clos}\,\left\{\gamma_\alpha^\Psi(\xi); \,\,\xi\in
G_2\right\}.$$ Moreover, for a real-valued function
$a(r,s)=\alpha(r)$ we have
$$\mathrm{sp}\,T_\alpha^\Psi = \left[\inf_{\xi\in G_2}\gamma_\alpha^\Psi(\xi),
\sup_{\xi\in G_2} \gamma_\alpha^\Psi(\xi)\right].$$
\end{thm}

Introduce the $C^*$-algebra $\mathcal{A}_\infty$ of all
$L_\infty(G, d\zeta)$-functions depending on $r$ only, where
$(r,s)\in G$. Then as a consequence of Theorem~\ref{CTO1} we have

\begin{thm}\label{coralgebra}
The $C^*$-algebra $\mathcal{T}_\Psi(\mathcal{A}_\infty)$ of TLO's
$T_a^\Psi$ with symbols $a\in \mathcal{A}_\infty$ is commutative
and is isometrically imbedded to the algebra $C_b(G_2)$ of bounded
continuous functions on $G_2$. The isomorphic imbedding
$$\tau^\Psi_\infty:
\mathcal{T}_\Psi(\mathcal{A}_\infty)\longrightarrow C_b(G_2)$$ is
generated by the following mapping $$\tau^\Psi_\infty: T_a^\Psi
\longmapsto \gamma_{a}^\Psi(\xi)$$ of generators of the algebra
$\mathcal{T}_\Psi(\mathcal{A}_\infty)$.
\end{thm}

Commutativity of the algebra
$\mathcal{T}_\Psi(\mathcal{A}_\infty)$ is a rather interesting
feature, see the book~\cite{vasilevskibook} devoted to this
phenomena for Toeplitz operators on the Bergman spaces. For two
symbols $a,b\in \mathcal{A}_\infty$ we have obviously that, in
general,
$$\gamma_{a}^\Psi(\xi)\gamma_{b}^\Psi(\xi)-\gamma_{ab}^\Psi(\xi)\neq
0,$$ which means that the $C^*$-algebra
$\mathcal{T}_\Psi(\mathcal{A}_\infty)$ generated by TLO's with
such symbols is a further example of algebra with the property
that for each pair $a,b\in \mathcal{A}_\infty$ the commutator
$[T_a^\Psi, T_b^\Psi]=0$, while the semi-commutator $[T_a^\Psi,
T_b^\Psi)$ is not compact.

Fix a number $m\in\mathbb{N}$. Let $Y_k^\Psi$, $k=1,\dots, m$, be
disjoint measurable sets in $G_1$ with a positive measure, such
that $\bigcup_{k=1}^m Y_k^\Psi=G_1$. Let $\Pi_k^\Psi =
G_2+\mathrm{i}\,Y_k^\Psi$, $k=1,\dots,m$, be the corresponding
sets in $G=G_2+\mathrm{i}\,G_1$. Denote by $\chi_{Y_k^\Psi}$ the
characteristic function of the set $Y_k^\Psi$, and by
$\chi_{\Pi_k^\Psi}$ the characteristic function of the set
$\Pi_k^\Psi$, $k=1,\dots, m$, respectively. For the algebra
$$\mathcal{A}_m^\Psi = \left\{a_1\chi_{\Pi_1^\Psi}+\dots+a_m\chi_{\Pi_m^\Psi}; \,\,
a_k\in\mathbb{C}, k=1,\dots, m\right\}$$ we immediately have the
following result.

\begin{thm}\label{thm3.6}
The algebra $\mathcal{T}_\Psi(\mathcal{A}_m^\Psi)$ is isometric
and isomorphic to the algebra $C(\nabla_\Psi)$, where
$$\nabla_\Psi = \nabla_\Psi\left(Y_1^\Psi, \dots, Y_m^\Psi\right) =
\mathrm{clos}\,\left\{\left(\gamma^\Psi_{\chi_{Y_1^\Psi}}(\xi),
\dots, \gamma^\Psi_{\chi_{Y_m^\Psi}}(\xi)\right); \,\,\xi\in
G_2\right\},$$ and the functions
$\gamma^\Psi_{\chi_{Y_k^\Psi}}(\xi)$, $k=1,\dots,m$, are given by
$$\gamma^\Psi_{\chi_{Y_k^\Psi}}(\xi) = \int_{G_1} \chi_{Y_k^\Psi}(z)
|\ell_\Psi(z,\xi)|^2\,\mathrm{d}\zeta_1(z) = \int_{Y_k^\Psi}
|\ell_\Psi(z,\xi)|^2\,\mathrm{d}\zeta_1(z), \quad \xi\in G_2.$$
The isomorphism
$$\tau_m^\Psi: \mathcal{T}_\Psi(\mathcal{A}_m^\Psi)\longrightarrow C(\nabla_\Psi)$$ is
generated by the following mapping of generators $T_a^\Psi$ of the
algebra $\mathcal{T}_\Psi(\mathcal{A}_m^\Psi)$
$$\tau_m^\Psi: T_a^\Psi \longmapsto a_1 z_1^\Psi+\dots+a_m z_m^\Psi, \quad
z^\Psi=\left(z_1^\Psi,\dots, z_m^\Psi\right)\in\nabla_\Psi,$$
where
$$a=a_1\chi_{\Pi_1^\Psi}+\dots+a_m\chi_{\Pi_m^\Psi}\in
\mathcal{A}_m^\Psi.$$
\end{thm}

Property to be unitarily equivalent to a multiplication operator
permits us to describe easily a sufficiently rich structure of
invariant subspaces of $C^*$-algebra
$\mathcal{T}_\Psi(\mathcal{A}_\infty)$.

\begin{thm}
The commutative $C^*$-algebra
$\mathcal{T}_\Psi(\mathcal{A}_\infty)$ is reducible and every
invariant subspace $\mathcal{S}_\Psi$ of
$\mathcal{T}_\Psi(\mathcal{A}_\infty)$ is defined by a measurable
subset $S_\Psi\subset G_2$ and has the form
$$\mathcal{S}_\Psi=(R_\Psi^*\chi_{S_\Psi} I)L_2(G_2, d\zeta_2)$$
with $\chi_{S_\Psi}$ being the characteristic function of
$S_\Psi$.
\end{thm}

\subsection{TLO's with symbols depending on the second
variable}\label{section3.2}

Now we are interested in symbols depending only on second
(vertical) variable of $G$. In this case the TLO is no more
unitarily equivalent to a multiplication operator, but certain
class of integral operators appears.

\begin{thm}\label{CTO2}
If a measurable function $a(r,s)=\beta(s)$ does not depend on $r$,
then the TLO $T_{\beta}^{\Psi}$ acting on
$W_\Psi(L_2(\mathbb{R}))$ is unitarily equivalent to the following
integral operator
$$\left(\mathfrak{B}_\beta^{\Psi}f\right)(\xi) =
\int_{G_2} \mathfrak{b}_\Psi(\xi,\omega)
\,\hat{\beta}(\pm(\xi-\omega))f(\omega)\,\mathrm{d}\zeta_2(\omega),
\quad \xi\in G_2,
$$ acting on $L_2(G_2, \mathrm{d}\zeta_2)$. The function $\mathfrak{b}_\Psi: G_2\times G_2
\to\mathbb{C}$ is given by
$$\mathfrak{b}_\Psi(\xi,\omega) = \int_{G_1}
\ell_\Psi(r,\omega)\overline{\ell_\Psi(r,\xi)}\,\mathrm{d}\zeta_1(r).$$
\end{thm}

\begin{proof} Similarly as in the proof of Theorem~\ref{CTO1} we get
$$\mathfrak{B}_\beta^{\Psi} = R_\Psi T^{\Psi}_{\beta(s)} R_\Psi^* = R_\Psi \beta(s) R_\Psi^* = Q_\Psi^* U_\Psi \beta(s)
U_\Psi^* Q_\Psi.$$ Using the convolution theorem for Fourier
transform we have {\setlength\arraycolsep{2pt}
\begin{align*}
\left(\mathfrak{B}_\beta^{\Psi} f\right)(\xi) & = Q_\Psi^* \left(
\int_{G_2} \hat{\beta}(\pm(s-\omega)) (Q_\Psi f)(r,\omega)\,\mathrm{d}\zeta_2(\omega)\right)(\xi) \\
& = \int_{G_1} \overline{\ell_\Psi(r,\xi)}\,\mathrm{d}\zeta_1(r)
\int_{G_2} \hat{\beta}(\pm(\xi-\omega)) f(\omega)\,\ell_\Psi(r,\omega)\,\mathrm{d}\zeta_2(\omega) \\
& = \int_{G_2} \hat{\beta}(\pm(\xi-\omega))
f(\omega)\,\mathrm{d}\zeta_2(\omega)
\int_{G_1} \ell_\Psi(r,\omega)\overline{\ell_\Psi(r,\xi)}\,\mathrm{d}\zeta_1(r) \\
& = \int_{G_2} \mathfrak{b}_\Psi(\xi,\omega)
\hat{\beta}(\pm(\xi-\omega)) f(\omega)
\,\mathrm{d}\zeta_2(\omega), \quad \xi\in G_2,
\end{align*}}where $f\in L_2(G_2, \mathrm{d}\zeta_2)$, and $$
\mathfrak{b}_\Psi(\xi,\omega) = \int_{G_1}
\ell_\Psi(r,\omega)\overline{\ell_\Psi(r,\xi)}\,\mathrm{d}\zeta_1(r),$$
which completes the proof. \end{proof}

\begin{rem}\rm
Observe that the class of integral operators
$$\left(\mathfrak{B}_\beta^{\Psi}f\right)(\xi) = \int_{G_2}
\mathfrak{b}_\Psi(\xi,\omega)
\hat{\beta}(\pm(\xi-\omega))f(\omega)\,\mathrm{d}\omega, \quad
\xi\in G_2,
$$ is interesting itself, and in some sense extends and generalizes the class of
operators considered in~\cite{cordes} to the whole line. Note that
for each $\Psi\in L_2(\mathbb{R})$ the function
$\mathfrak{b}_\Psi$ has the following properties:
$$\mathfrak{b}_\Psi(\xi,\omega)=\overline{\mathfrak{b}_\Psi(\omega,\xi)}\,\,\,\textrm{for
all}\,\,\, \xi, \omega\in G_2, \qquad
\mathfrak{b}_\Psi(\xi,\xi)=1\,\,\,\textrm{for a.e.}\,\,\,\xi\in
G_2.$$ Further remarkable properties may be obtained when
considering some special cases of wavelets, or windows,
respectively.
\end{rem}

\subsection{General case of symbols}\label{section3.3}

The above mentioned construction of unitary operators may be used
to study more general symbols $a(r,s)$ for which the TLO
$T_a^\Psi$ is no longer unitarily equivalent to a multiplication
operator, because the operator $R_\Psi T_a^\Psi R_\Psi^*$ might
have a more complicated structure as we have demonstrated above
for the case of symbols depending on the second variable. As we
will prove now the TLO $T_a^\Psi$ with symbols which depend on
both variables $(r,s)\in G$ is unitarily equivalent to a
pseudodifferential operator with certain compound (double) symbol.
We clarify this statement for the case of symbol $a$ in the
product form $a(r,s)=\alpha(r)\beta(s)$.

\begin{thm}\label{CTO3}
Let $a(r,s)=\alpha(r)\beta(s)$ be a measurable symbol on $G$. Then
the TLO $T_{a}^{\Psi}$ acting on $W_\Psi(L_2(\mathbb{R}))$ is
unitarily equivalent to the pseudodifferential operator
$\mathfrak{A}^\Psi_{\mathfrak{a}}$ acting on $L_2(G_2,
\mathrm{d}\zeta_2)$. The operator
$\mathfrak{A}^\Psi_{\mathfrak{a}}$ is given by the iterated
integral
\begin{equation}\label{pseudodiff}
\left(\mathfrak{A}^\Psi_{\mathfrak{a}} f\right)(x) = \int_{G_2}
\mathrm{d}\zeta_2(\xi) \int_{G_2} \mathfrak{a}_\Psi(x,y,\xi)
f(y)\,\mathrm{e}^{\mp 2\pi
\mathrm{i}(x-y)\xi}\,\mathrm{d}\zeta_2(y),
\end{equation} for $x\in G_2$, where the compound (double) symbol $\mathfrak{a}_\Psi:
G_2\times G_2\times G_2\to\mathbb{C}$ has the form
$\mathfrak{a}_\Psi(x,y,\xi)= \Gamma^\Psi_\alpha(x,y) \beta(\xi)$
with
$$\Gamma^\Psi_\alpha(x,y) = \int_{G_1} \alpha(r)
\overline{\ell_\Psi(r,x)}\ell_\Psi(r,y)\,\mathrm{d}\zeta_1(r).$$
\end{thm}

Similarly as above, since we consider the operator
$U_\Psi=\left(I\otimes \mathcal{F}^{\pm 1}\right)$ to describe the
both cases of wavelet and time-frequency analysis, also here the
signs $\mp$ in~(\ref{pseudodiff}) correspond to
$\mathfrak{A}^\psi_{\mathfrak{a}}$ (for wavelet case) and
$\mathfrak{A}^\phi_{\mathfrak{a}}$ (for time-frequency case),
respectively.

\begin{proof} Let $f\in L_2(G_2, \mathrm{d}\zeta_2)$. Then $T_a^\Psi$ is unitarily
equivalent to the operator $$\mathfrak{A}_{\mathfrak{a}}^{\Psi} =
R_\Psi T^{\Psi}_{\alpha(r)\beta(s)} R_\Psi^* = R_\Psi
\alpha(r)\beta(s) R_\Psi^* = Q_\Psi^* \alpha(r) U_\Psi \beta(s)
U_\Psi^* Q_\Psi.$$ Then
\begin{align*}
(\mathfrak{A}^\Psi_{\mathfrak{a}} f)(\lambda) & = \int_{G_1}
\alpha(r) \overline{\ell_\Psi(r,\lambda)}\,\mathrm{d}\zeta_1(r)
\int_{G_2} \beta(s)\,\mathrm{e}^{\mp 2\pi
\mathrm{i}s\lambda}\,\mathrm{d}\zeta_2(s)
\\ & \phantom{=} \times \int_{G_2} f(\omega)\ell_\Psi(r,\omega)\,
\mathrm{e}^{\pm 2\pi \mathrm{i}\omega s}\,\mathrm{d}\zeta_2(\omega) \\
& = \int_{G_2} \beta(s)\,\mathrm{e}^{\mp 2\pi
\mathrm{i}s\lambda}\,\mathrm{d}\zeta_2(s) \int_{G_2}
f(\omega)\,\mathrm{e}^{\pm 2\pi \mathrm{i}\omega
s}\,\mathrm{d}\zeta_2(\omega)
\\ & \phantom{=} \times \int_{G_1} \alpha(r)\overline{\ell_\Psi(r,\lambda)}
\ell_\Psi(r,\omega)\,\mathrm{d}\zeta_1(r).
\end{align*}If the last integral is denoted by $\Gamma^\Psi_\alpha(\lambda,\omega)$,
then we finally have
\begin{align*}
\left(\mathfrak{A}_{\mathfrak{a}}^\Psi f\right)(\lambda) & =
\int_{G_2} \mathrm{d}\zeta_2(s) \int_{G_2}
\Gamma^\Psi_\alpha(\lambda,\omega)\beta(s) f(\omega)\,
\mathrm{e}^{\mp 2\pi
\mathrm{i}(\lambda-\omega)s}\,\mathrm{d}\zeta_2(\omega)
\\ & = \int_{G_2} \mathrm{d}\zeta_2(s) \int_{G_2}
\mathfrak{a}_\Psi(\lambda,\omega,s) f(\omega)\,\mathrm{e}^{\mp
2\pi \mathrm{i}(\lambda-\omega)s}\,\mathrm{d}\zeta_2(\omega).
\end{align*} Changing the variables $\lambda=x$, $\omega=y$ and $s=\xi$
we finally get the standard notation~(\ref{pseudodiff}) for the
pseudodifferential operator with a compound symbol, see
e.g.~\cite{karlovich}.
\end{proof}

\begin{rem}\rm
Observe that $\Gamma^\Psi_\alpha(x,x) = \gamma^\Psi_\alpha(x)$ is
the function from Theorem~\ref{CTO1} which is responsible for
properties of the TLO $T_\alpha^\Psi$ whose symbol
$a(r,s)=\alpha(r)$ depends on the first variable only. Also, for
$\alpha(r)\equiv 1$ on $G_1$ we get $\Gamma^\Psi_1(x,y) =
\mathfrak{b}_\Psi(x,y)$, the function appearing in
Theorem~\ref{CTO2}. Note that in the proof of Theorem~\ref{CTO3}
we did not use the convolution theorem as in Theorem~\ref{CTO2} to
get the desired form~(\ref{pseudodiff}) of pseudodifferential
operator. Thus, the operator $\mathfrak{B}_\beta^\Psi$ may also be
viewed as a pseudodifferential operator of the
form~(\ref{pseudodiff}). Further research on TLO's $T_a^\Psi$
using the deeper connection with pseudodifferential operators will
be considered elsewhere.
\end{rem}

\section*{Concluding remarks}

In the end let us note that the presented technique and the
obtained results are interesting from various viewpoints, because
(inter alia)

(i) they represent a unified approach to study both
Calder\'on-Toeplitz and Gabor-Toeplitz operators and properties
which they share;

(ii) the presented technique is purely analytic based only on
operator theory and does not use neither the specifics of groups
behind the construction of localization operators (affine, or
Weyl-Heisenberg group, respectively), nor time-scale, or
time-frequency methods;

(iii) they give rise to commutative algebras of TLO's which are
practically unknown;

(iv) they enable to study the TLO's using their unitarily
equivalent images in the class of pseudodifferential operators
with compound symbols as an analog of the Berezin approach known
for Toeplitz and Weyl pseudodifferential operators.

\vspace{0.5cm}{\bf Acknowledgements:} The first author has been on
a postdoctoral stay at the Departamento de Matem\'aticas,
CINVESTAV del IPN (M\'exico), when writing this paper and
investigating the topics presented herein. He therefore gratefully
acknowledges the hospitality and support of the mathematics
department of CINVESTAV on this occasion. Especially, he would
like to thank Professor Nikolai L. Vasilevski for his
encouragement and help during his stay in M\'exico.

\vspace{5mm}

\noindent \small{Ondrej Hutn\'ik, Departamento de Matem\'aticas,
CINVESTAV del IPN, {\it Current address:} Apartado Postal 14-740,
07000, M\'exico, D.F., M\'exico
\newline {\it E-mail address:} hutnik@math.cinvestav.mx}
\newline AND \newline \noindent \small{Institute of Mathematics,
Faculty of Science, Pavol Jozef \v Saf\'arik University in Ko\v
sice, Jesenn\'a 5, 040~01 Ko\v sice, Slovakia
\newline {\it E-mail address:} ondrej.hutnik@upjs.sk}

\vspace{5mm} \noindent \small{M\'aria Hutn\'ikov\'a, Department of
Physics, Faculty of Electrical Engineering and Informatics,
Technical University of Ko\v sice, {\it Current address:} Park
Komensk\'eho~2, 042~00 Ko\v sice, Slovakia \newline {\it E-mail
address:} maria.hutnikova@tuke.sk}

\end{document}